\documentclass[11pt]{article}
\usepackage{latexsym}
\usepackage{amsfonts}
\usepackage{graphicx}
\usepackage[sort]{natbib}
\usepackage{lineno}
\usepackage{amsmath}
\usepackage{amssymb}
\usepackage{xcolor}
\usepackage{cancel}

\newcommand{\coln}{\hspace*{-6pt}{\bf :}}
\newtheorem{theorem}{Theorem}
\newtheorem{lemma}{Lemma}
\newtheorem{corol}{Corollary}
\newtheorem{property}{Property}
\newtheorem{conjecture}{Conjecture}

\newcommand{\Theorem}[1]{\begin{theorem}\coln ~#1 \end{theorem}}

\newcommand{\Corol}[1]{\begin{corol}\coln ~#1 \end{corol}}
\newcommand{\proof}[1]{\noindent{\bf Proof:~}#1\hfill ~$\Box$}
\newcommand{\be}{\begin{equation}}
\newcommand{\ee}{\end{equation}}

\newcommand{\Obar}{\overline{\Omega}}

\begin{document}
\title{\bf On the Combined Inverse-Square Effect of Multiple Point Sources in Multidimensional Space}
	\author{
		Keaton Coletti\\College of Engineering\\ University of Georgia - Athens\\ Athens, GA 30602\\ e-mail: keaton.coletti@uga.edu
		\and Pawel Kalczynski\\
		College of Business and Economics\\
		California State University - Fullerton\\
		Fullerton, CA 92834.\\e-mail: pkalczynski@fullerton.edu\and	Zvi Drezner\\
		College of Business and Economics\\
		California State University - Fullerton\\
		Fullerton, CA 92834.\\e-mail: zdrezner@fullerton.edu
	}
\date{}
\maketitle
\begin{abstract}
The inverse-square law states that the effect a source has on its surroundings is inversely proportional to the square of the Euclidean distance from that source. Its applicability spans multiple fields including physics, engineering, and computer science. We study the combined effect of multiple point sources surrounding a closed region in multidimensional space. We determine that the maximum effect in D dimensions is on the region's boundary if $D \leq 4$, and the minimum is on the boundary if $D \geq 4$.
\end{abstract}
\renewcommand{\baselinestretch}{1.5}
\renewcommand{\arraystretch}{0.66}
\large
\normalsize

\section{Introduction}
This paper investigates Euclidean distances in multidimensional space. Point sources emit some quantity whose effect is inversely proportional to the squared distance from the source, and the total effect at a point in space is the sum of scalar contributions from each source. This scalar sum applies to radiation, sound, and light.
The cooperative cover location problem is proposed in \citet{BDK09}. There are sources of light, cell phone signals, noise, heat, etc. that decline by the square of the distance. Each point in a region is affected by the sum of the signals. For example, light posts in a parking lot are in three-dimensional space, and we are interested in the total light in the two-dimensional parking lot. Which point in the parking lot receives the minimum amount of light? This minimum value guarantees that all the points in the parking lot receive at least this minimum light. 

Locating global extrema of such quantities is often of interest \citep{HPT81,DDK20}. The results presented in this work reduce optimization dimensionality by 1, thereby improving computation time and reducing the risk of erroneously identifying local extrema.

As an example, consider one or more closed shapes in 3 dimensions surrounded by point sources. The sources emit sound, heat, light, or other types of energy in every direction. The emission from these sources may have different levels of intensity. If the optimization objective is to arrange the sources and/or the shapes in a manner that reduces the maximum total (cumulative) effect of the emissions on the shapes (keeping all sources outside of the shapes), we show that it is sufficient to consider solely the effect on the boundaries of these shapes because, in 3 dimensional space, the maximum total effect must be on the boundary.

\section{Problem Statement}
    Let $D\in\{1,2,3,...\}$ be the dimension of a Euclidean space $\mathbb{R}^D$. Consider $n$ point sources $S_i=(s_{1i},s_{2i},...,s_{Di})\in \mathbb{R}^D$, $i=1,2,...,n$. Define $w_j>0$ which represent the strength of the quantity emitted by source $i$.
    
    Furthermore, let $\Obar\in\mathbb{R}^D$, be a region affected by the point sources, such that $\Obar\, \cap\, \{S_1,S_2,...,S_n\}=\varnothing$. The combined effect of point sources on $X=(x_1,x_2,...,x_D)$ is given by
    
    \begin{equation}
    J(X)=\sum\limits_{i=1}^n\frac{w_i}{\sum\limits_{k=1}^D (x_k-s_{ki})^2}. \label{obj}    
    \end{equation}

\section{Maximum and Minimum Values}

\Theorem{\label{th1}
Let $\Obar \in \mathbb{R}^D\setminus \{S_i\}$ contain the open set $\Omega$ and its boundary $\partial \Omega$. If $J(X)$ is non-constant on $\Omega$,
\begin{enumerate}
    \item If $D \leq 4$, then $J(X)$ obtains maximum in $\Obar$ only on $\partial \Omega$.
    \item If $D \geq 4$, then $J(X)$ obtains minimum in $\Obar$ only on $\partial \Omega$.
\end{enumerate}
}

\proof{
$J(X)$ is defined on $\Obar$ and hence attains a maximum and minimum in $\Obar$. $J(X)$ is a continuous function and has continuous derivatives of every order in $\Obar$. To simplify notation, define $x_{ki}=(x_k-s_{ki})$ for $k=1,2,...,D, i=1,2,...,n$ (note that $dx_{ki}=dx_k$). Then Eq. \eqref{obj} can be written as

$$
J(X)=\sum\limits_{i=1}^n\frac{w_i}{\sum\limits_{k=1}^Dx_{ki}^2},
$$

and its unmixed first and second-partial derivatives are

$$
\frac{dJ(X)}{dx_k}=\sum\limits_{i=1}^n\frac{-2x_{ki}w_i}{\left[\sum\limits_{k=1}^Dx_{ki}^2\right]^2}
$$

$$
\frac{d^2J(X)}{dx_k^2}=\sum\limits_{i=1}^n\frac{-2w_i\left[\sum\limits_{k=1}^Dx_{ki}^2\right]^2+8x_{ki}^2w_i\sum\limits_{k=1}^Dx_{ki}^2}{\left[\sum\limits_{k=1}^Dx_{ki}^2\right]^4}.
$$

Computing the Laplacian yields

\begin{equation}
\Delta J(X) = \sum\limits_{k=1}^D\frac{d^2J(X)}{dx_k^2}=\sum\limits_{i=1}^n\frac{-2Dw_i\left[\sum\limits_{k=1}^Dx_{ki}^2\right]^2+8w_i\left[\sum\limits_{k=1}^Dx_{ki}^2\right]^2}{\left[\sum\limits_{k=1}^Dx_{ki}^2\right]^4}=\sum\limits_{i=1}^n\frac{(8-2D)w_i}{\left[\sum\limits_{k=1}^Dx_{ki}^2\right]^2}. \label{laplace}
\end{equation}

We have two cases:
\begin{enumerate}
    \item $\Delta J(X) \geq 0$ for $D \leq 4$, so $J$ is sub-harmonic For all $X \in \Obar$, $\frac{d^2J(X)}{dx_k^2} > 0$ for some $k$, so $J$ has no local maxima or plateaus in $\Obar$.
    \item $\Delta J(X) \leq 0$ for $D \geq 4$, so $J$ is super-harmonic For all $X \in \Obar$, $\frac{d^2J(X)}{dx_k^2} < 0$ for some $k$, so $J$ has no local minima or plateaus in $\Obar$.
\end{enumerate}
The maximum principle for sub-harmonic functions (case 1) and minimum principle for super-harmonic functions (case 2) is applied by \citet{stein,axler,evans,ovall}.
}

\Corol{\label{cor1}
If $D=4$, $J$ is harmonic on its domain and attains both maximum and minimum values in $\Obar$ on the boundary.
}

Note that $\Obar$ need not be a connected domain. No statement has yet been made about the maximum of $J(X)$ for $D>4$ or the minimum of $J(X)$ for $D<4$. Theorem \ref{th1} does not apply for other values of $D$.

\subsection{When Maximum and Minimum Values Can Be in the Interior}
We construct a counter example for which in $D\ge 5$ the maximum can be in the interior of $\Obar$, and in $D\le 3$ the minimum can be in the interior.  Consider a set of $2D$ source points in D-dimensional space, as depicted in Table \ref{counter} for $D=5$, with $w_i=1$, $i=1,2,..,2D$. All the points are on the surface of a D-dimensional sphere centered at the origin (0, 0, $\ldots$, 0) with a radius of $\frac12$. We first show that the center is a local maximum for $D\ge 5$, and a local minimum for $D\le 3$. We then find a radius $0<r_{max} <\frac12$ of a D-dimensional sphere that defines the region $\Obar$.

\begin{table}
	\begin{center}
		\caption{\label{counter}Counter example for $D=5$}
		\begin{tabular}{|c|c|c|c|c|c|}
			\hline
			&$k=1$&$k=2$&$k=3$&$k=4$&$k=5$\\\hline
			$i=1$&$\frac12$&0&0&0&0\\\hline
			$i=2$&0&$\frac12$&0&0&0\\\hline
			$i=3$&0&0&$\frac12$&0&0\\\hline
			$i=4$&0&0&0&$\frac12$&0\\\hline
			$i=5$&0&0&0&0&$\frac12$\\\hline
			$i=6$&$-\frac12$&0&0&0&0\\\hline
			$i=7$&0&$-\frac12$&0&0&0\\\hline
			$i=8$&0&0&$-\frac12$&0&0\\\hline
			$i=9$&0&0&0&$-\frac12$&0\\\hline
			$i=10$&0&0&0&0&$-\frac12$\\\hline
		\end{tabular}
	\end{center}
\end{table}

Consider a location $\Delta X=\{\Delta x_k\}$ for $k=1,\ldots, D$ on a D-dimensional sphere of radius $r$ centered at the origin so that $\sum\limits_{k=1}^D\Delta x_k^2 = r^2$. For source point $i=1$,
$\sum\limits_{k=1}^D(\Delta x_k-x_i)^2=(\frac12-\Delta x_1)^2+\sum\limits_{k=2}^D\Delta x_k^2=\frac14-\Delta x_1+r^2$. Similarly, for $i=D+1$, $\sum\limits_{k=1}^D(\Delta x_k-x_i)^2=(-\frac12-\Delta x_1)^2+\sum\limits_{k=2}^D\Delta x_k^2=\frac14+\Delta x_1+r^2$, and similar expressions for the other points. Therefore,

\begin{eqnarray}
J(\Delta X)&=&\sum\limits_{i=1}^D\left[\frac{1}{\frac14-\Delta x_i+r^2}+\frac{1}{\frac14+\Delta x_i+r^2}\right]=\sum\limits_{i=1}^D\frac{\frac12+2r^2}{(\frac14+\Delta x_i+r^2)(\frac14-\Delta x_i+r^2)}\nonumber\\
&=&\sum\limits_{i=1}^D\frac{\frac12+2r^2}{(\frac14+r^2)^2-\Delta x_i^2}~.\nonumber
\end{eqnarray}
Since $\Delta x_i\le r$, we can get an upper bound for $J(\Delta X)$ by replacing $\Delta x_i$ by $r$, and a lower bound by replacing $\Delta x_i$ by 0. However, the gap in this approximation is too large and has an error of the order $r^2$ which fails when used to prove local optima. A better approximation of order $r^4$ is obtained by multiplying both numerator and denominator by $(\frac14+r^2)^2+\Delta x_i^2$ getting:
\begin{equation}\label{JDX}
J(\Delta X)=2\sum\limits_{i=1}^D\frac{\left[\frac14+r^2\right]\left[(\frac14+r^2)^2+\Delta x_i^2\right]}{(\frac14+r^2)^4-\Delta x_i^4}
\end{equation}
and getting bounds by replacing $\Delta x_i$ by $r$ or 0 {\it only in the denominator}, and applying the property that $\sum\limits_{i=1}^D\Delta x_i^2=r^2$.

\subsection*{Upper Bound:}
\begin{eqnarray}
J(\Delta X)&\le&2\sum\limits_{i=1}^D\frac{\left[\frac14+r^2\right]\left[(\frac14+r^2)^2+\Delta x_i^2\right]}{(\frac14+r^2)^4-r^4}\nonumber\\&=&
2D\frac{\left[\frac14+r^2\right](\frac14+r^2)^2}{(\frac14+r^2)^4-r^4}+
2\sum\limits_{i=1}^D\frac{\left[\frac14+r^2\right]\Delta x_i^2}{(\frac14+r^2)^4-r^4}
=2D\frac{(\frac14+r^2)^3}{(\frac14+r^2)^4-r^4}+
2\frac{\left[\frac14+r^2\right]r^2}{(\frac14+r^2)^4-r^4}\nonumber\\&=&\frac{\frac14+r^2}{(\frac14+r^2)^4-r^4}\left[\frac{D}8+(D+2)r^2+2Dr^4\right]~.\label{UB}
\end{eqnarray}
\subsection*{Lower Bound:}
\begin{eqnarray}
J(\Delta X)&\ge&2\sum\limits_{i=1}^D\frac{\left[\frac14+r^2\right]\left[(\frac14+r^2)^2+\Delta x_i^2\right]}{(\frac14+r^2)^4}\nonumber\\
&=&\frac{2D}{\frac14+r^2}+
2\sum\limits_{i=1}^D\frac{\Delta x_i^2}{(\frac14+r^2)^3}
=\frac{2D}{\frac14+r^2}+
\frac{2r^2}{(\frac14+r^2)^3}\nonumber\\&=&\frac{1}{(\frac14+r^2)^3}\left[\frac{D}8+(D+2)r^2+2Dr^4\right]~.\label{LB}
\end{eqnarray}

\Theorem{\label{th2}
    If $D > 4$ (resp. $D < 4$), $J(X)$ can attain maximum (resp. minimum) in $\Obar$ on its interior or boundary.
}
\proof{ To show that the optimum can be in the interior of  $\Obar$, consider $\Obar$ as a D-dimensional sphere of a small radius $r$ centered at the origin.
	
Since $J(0)=8D$, by equation (\ref{UB})
\begin{eqnarray}
J(\Delta X)-J(0)&\le&\frac{\frac14+r^2}{(\frac14+r^2)^4-r^4}\left[\frac{D}8+(D+2)r^2+2Dr^4\right]-8D\nonumber\\
&=&\frac{(\frac14+r^2)\left[\frac{D}8+(D+2)r^2+2Dr^4\right]-8D((\frac14+r^2)^4-r^4)}{(\frac14+r^2)^4-r^4}\nonumber\\
&=&\frac{-8Dr^8-6Dr^6+\frac{13D+4}{2}r^4+\frac{4-D}{8}r^2}{(\frac14+r^2)^4-r^4}~,\label{DiffUB}
\end{eqnarray}
and by equation (\ref{LB})
\begin{eqnarray}
J(\Delta X)-J(0)&\ge&\frac{1}{(\frac14+r^2)^3}\left[\frac{D}8+(D+2)r^2+2Dr^4\right]-8D\nonumber\\
&=&\frac{\frac{D}8+(D+2)r^2+2Dr^4-8D(\frac14+r^2)^3}{(\frac14+r^2)^3}\nonumber\\
&=&\frac{-8Dr^6-4Dr^4+\frac{4-D}{2}r^2}{(\frac14+r^2)^3}~.\label{DiffLB}
\end{eqnarray}

To show that $J(\Delta X)< J(0)$ for $r>0$ in the interior of $\Obar$, it suffices to show that the polynomial at the numerator of (\ref{DiffUB}) is negative for $r>0$.
The first derivative of the polynomial at $r=0$ is 0, and the second derivative is negative for $D\ge 5$. Therefore, the polynomial is negative up to a certain value of $r$ for $D\ge 5$.
To show that $J(\Delta X)> J(0)$ for $r>0$ in the interior of $\Obar$, it suffices to show that the polynomial at the numerator of (\ref{DiffUB}) is positive for $r>0$.
The first derivative of the polynomial at $r=0$ is 0, and the second derivative is positive for $D\le 3$. Therefore, the polynomial is positive up to a certain value of $r$ for $D\le 3$.}

\subsection{A radius of assured decrease for $D \geq 5$, and increase for $D\leq 3$, in the counter example}
Consider the counter example in Theorem  \ref{th2}. The proof of Theorem \ref{th2} showed that for $D>4$, there exists a D-dimensional sphere of radius $r\leq r_{\max}$, for some $r_{\max}$, centered at the origin on which $J(X)$ has guaranteed global maximum in $\Obar$ at the origin, and for $D<4$ there exist a global minimum.

To find a lower bound for the radius of the D-dimensional sphere for $D\ge 5$, we find the smallest positive root $\bar x$ of the polynomial on the numerator of equation (\ref{DiffUB}) divided by $r^2$, where $x=r^2$:

$$
-8Dx^3-6Dx^2+\frac{13D+4}{2}x+\frac{4-D}{8}~=~0~,
$$
and $r_{\max} =\sqrt{\bar x}$. By numerical solution we found lower bounds for $r_{\max}$ values depicted in Table~\ref{r}. A good approximation (which is a bit lower) is obtained by ignoring $x^2$ and $x^3$ that have negative coefficients. This yields: $r_{\max}\approx \sqrt{\frac{D-4}{52D+16}}$ which converges to $\sqrt{\frac{1}{52}}=0.1387$ as $D\rightarrow\infty$. For $D=1,000,000$ the root of the polynomial yields $r_{\max}=0.1400$. If we ignore just the term $x^3$, we get by solving the quadratic equation $r_{\max}=\sqrt{\frac{13D+4-\sqrt{157D^2+152D+16}}{24D}}$ which agrees with the values in Table \ref{r} to four significant digits. The limit as $D\rightarrow\infty$ is $\sqrt{\frac{13-\sqrt{157}}{24}}=\frac{1}{\sqrt{26+\sqrt{628}}}$.

A similar estimate for $r_{\max}$ for the minimum case can be developed by the polynomial in the numerator of (\ref{DiffLB}) yielding $r_{\max}=\sqrt{\frac{1}{\sqrt{4D}}-\frac14}$. For $D=2$, $r_{\max}=0.3218$, and for $D=3$ $r_{\max}=0.1967$. Note that all these values, for both cases, are less than 0.5 so that all the source points are outside  $\Obar$.

\begin{table}
	\begin{center}
		\caption{\label{r}Limit on maximum $r$}
		\begin{tabular}{|c|c||c|c|}
			\hline
			$D$&$r_{\max}$&$D$&$r_{\max}$\\\hline
5&	0.0603&	10&	0.1063\\
6&	0.0783&	15&	0.1183\\
7&	0.0892&	20&	0.1240\\
8&	0.0966&	50&	0.1337\\
9&	0.1021&	100&	0.1369\\
\hline
		\end{tabular}
	\end{center}
\end{table}

\section{Conclusions}

We considered the combined effect of multiple source points emitting some quantity according to the inverse-square law, on an open domain that does not include these points. We showed that the combined (total) effect of point sources on the domain enjoys the following properties:
\begin{enumerate}
    \item The maximum must be on the boundary (not inside) for $D \leq 4$
    \item The minimum must be on the boundary (not inside) for $D \geq 4$
    \item The maximum may be inside the domain for $D > 4$
    \item The minimum may be inside the domain for $D < 4$.
\end{enumerate}

Our findings enable dimensionality reduction in practical optimization problems, resulting in significant efficiency gains and robustness.

\renewcommand{\baselinestretch}{1}
\renewcommand{\arraystretch}{1}
\large
\normalsize

\bibliographystyle{APAlike}


\end{document}